\documentclass{scrartcl}
\usepackage[numbers]{natbib}
\usepackage{fontenc}
\usepackage{shadethm}
\usepackage{amsthm}
\usepackage{float}
\usepackage{bm}
\usepackage{stmaryrd}
\usepackage{mathrsfs} 
\usepackage{amsmath}
\usepackage{amssymb}
\usepackage{graphicx}
\usepackage{subfigure}
\usepackage{wasysym}
\usepackage{sectsty}
\usepackage{hyperref}
\usepackage[a4paper, left=2.7cm, right=2.7cm, top=2.5cm]{geometry}
\usepackage{chngcntr}
\counterwithin*{equation}{section}

\usepackage{cleveref}
\DeclareMathAlphabet{\mathpzc}{OT1}{pzc}{m}{it}

\sectionfont{\normalsize\centering}
\subsectionfont{\footnotesize\centering}

\theoremstyle{plain}
\newtheorem{thm}{Theorem}[section] 

\theoremstyle{definition}
\newtheorem{defn}[thm]{Definition} 
\newtheorem{exmp}[thm]{Example} 
\newtheorem{lem}[thm]{Lemma}

\newtheorem{cor}[thm]{Corollary}

\def\XXint#1#2#3{{\setbox0=\hbox{$#1{#2#3}{\int}$ }
		\vcenter{\hbox{$#2#3$ }}\kern-.6\wd0}}

\usepackage[utf8]{inputenc}
\usepackage[german,english,russian]{babel}

\newcounter{MPequ}

\begin{document}\selectlanguage{english}
\begin{center}
\normalsize \textbf{\textsf{Zero set structure of real analytic Beltrami fields}}
\end{center}
\begin{center}
	Wadim Gerner\footnote{\textit{E-mail address:} \href{mailto:gerner@eddy.rwth-aachen.de}{gerner@eddy.rwth-aachen.de}}
\end{center}
\begin{center}
{\footnotesize	RWTH Aachen University, Lehrstuhl f\"ur Angewandte Analysis, Turmstra{\ss}e 46, D-52064 Aachen, Germany}
\end{center}
{\small \textbf{Abstract:} In this paper we prove a classification theorem for the zero sets of real analytic Beltrami fields. Namely, we show that the zero set of a real analytic Beltrami field on a real analytic, connected $3$-manifold without boundary is either empty after removing its isolated points or can be written as a countable, locally finite union of differentiably embedded, connected $1$-dimensional submanifolds with (possibly empty) boundary and tame knots.
	Further we consider the question of how complicated these tame knots can possibly be. To this end we prove that on the standard (open) solid toroidal annulus in $\mathbb{R}^3$, there exist for any pair $(p,q)$ of positive, coprime integers countable infinitely many distinct real analytic metrics such that for each such metric there exists a real analytic Beltrami field, corresponding to the eigenvalue $+1$ of the curl operator, whose zero set is precisely given by a standard $(p,q)$-torus knot. The metrics and the corresponding Beltrami fields are constructed explicitly and can be written down in Cartesian coordinates by means of elementary functions alone.}
\section{Introduction}
Beltrami fields on an oriented, real analytic Riemannian $3$-manifold $(M,g)$ without boundary are vector fields $\bm{X}$ which satisfy $\text{div}(\bm{X})=0$ and $\text{curl}(\bm{X})=\lambda \bm{X}$ for some smooth function $\lambda:M\rightarrow \mathbb{R}$. A special case of these are eigenvector fields of the curl operator corresponding to non-zero eigenvalues, since such eigenfields are automatically divergence-free. Such vector fields appear naturally in physics and have been widely studied in mathematics. For instance, they appear as stationary magnetic fields of the equations of ideal magnetohydrodynamics, and hence in particular in astrophysics, in the case of constant pressure and a resting plasma, \cite[Chapter III \S 1.A]{AK98}. On the other hand, they also appear as stationary solutions of the incompressible Euler equations for an appropriate pressure function, \cite[Chapter II \S 1.A]{AK98}. From a variational point of view, Beltrami fields (with a constant, non-zero proportionality function $\lambda$) are closely related to the helicity constraint magnetic energy minimisation, see \cite{W58}, \cite{A74}, \cite{AL91}, while Beltrami fields with a non-constant proportionality function were studied for instance in \cite{EP16} and \cite{N14}. In view of hydrodynamics Beltrami fields are of particular interest from a topological point of view. Namely, if we consider stationary solutions of the incompressible Euler equations, Arnold's theorem, \cite{A66}, \cite{A74}, \cite[Chapter II Theorem 1.2]{AK98}, in essence characterises the field line behaviour of any such solution unless it is everywhere collinear with its curl. Therefore if one is interested in 'complicated' field line behaviour of steady Euler flows one necessarily needs to consider Beltrami flows. For example, a consequence of Arnold's structure theorem is that if some real analytic nowhere vanishing, incompressible, steady Euler flow admits a 'chaotic' field line, i.e. a field line not contained in a codimension $1$ subset, then the flow is necessarily Beltrami, \cite[Chapter II, Proposition 6.2]{AK98}. That Beltrami flows indeed can have very interesting behaviour is well-known, see for example \cite{DFGHMS86}, \cite{EP12}, \cite{EP15}, \cite{EtGr00b}. The existence of 'knotted' field lines, at least for nowhere vanishing Beltrami fields with nowhere vanishing proportionality function, on closed $3$-dimensional manifolds is guaranteed by the hydrodynamical interpretation of the (proven) Weinstein conjecture, see \cite{EtGr00a} for the hydrodynamical interpretation and \cite{H093}, \cite{Ta07} for a proof of the Weinstein conjecture. In particular, this result tells us that each such Beltrami flow admits a closed field line, which is then necessarily a smoothly embedded circle, i.e. a tame knot. Observe that this result specifically assumes that the zero set of the underlying Beltrami field is empty. Contrary to the results regarding the (non-constant) field line behaviour of Beltrami flows, the structure of the zero set seems to have been investigated far less.  Let us assume for the moment that the proportionality function $\lambda$ is constant and that $M=\Omega$ is a domain in $\mathbb{R}^3$, then we in particular observe that
\[
-\Delta \bm{X}=\text{curl}(\text{curl}(\bm{X}))-\text{grad}(\text{div}(\bm{X}))=\lambda^2\bm{X},
\]
i.e., each component of the Beltrami field is in particular an eigenfunction of the Laplacian and consequently the zero set of $\bm{X}$ is the intersection of $3$ zero sets of eigenfunctions of $\Delta$. The zero sets, also referred to as nodal sets, of Laplacian eigenfunctions were thoroughly studied for example in \cite{HarSi89}, \cite{Log18}, \cite{SoZel11} and \cite{SoZel12}. These nodal sets are, modulo a codimension $2$ countably rectifiable subset, codimension $1$ hypersurfaces, \cite{HarSi89}. Thus, a priori, the zero set of a Beltrami field might be $2$-dimensional. However, it was for instance observed in \cite{C99}, that the zero sets of rotationally symmetric Beltrami fields on rotationally symmetric, bounded domains (diffeomorphic to the solid torus) are either empty (in case of the first eigenfield), \cite[Theorem 7]{C99}, or else they are 'well-separated' circles, \cite[Theorem 8]{C99}. In particular, the Hausdorff dimension is an integer not greater than $1$ and the zero sets, in the latter case, are circles, i.e. have a very special topological structure. In the present paper we will show that these features are more generally true for any real analytic Beltrami field defined on an abstract manifold without boundary. More precisely, our result states that, after removing the isolated points, the remaining nodal set is either empty or a countable, locally finite union of analytically embedded $1$-manifolds with (possibly non-present) $C^1$-endpoints, see \cref{MD1}, and tame knots. Our approach differs from the approach in \cite{C99}, since we do not assume any symmetry, but instead rely on results from semianalytic geometry, \cite{L65}, \cite{G68}, \cite{BM88}, most notably the curve selection lemma, \cite[\S 2]{Loo84}, \cite[Lemma 3.1]{M68}, \cite[Lemma 6.6]{MaTr07}.
\section{Main results}
\textbf{Conventions:} All manifolds are assumed to be Hausdorff, second countable, oriented, connected, real analytic and without boundary, unless otherwise noted. We will simply say: 'Let $M$ be a $3$-manifold', meaning it has all the previously mentioned properties and is $3$-dimensional. $\mathcal{V}^{\omega}(M)$ denotes the set of all real analytic vector fields on a given $3$-manifold $M$. Given a smooth Riemannian metric $g$ on a $3$-manifold $M$, we define the curl and divergence of a vector field by means of its identification with their corresponding $1$-form. More precisely if $\bm{X}$ is any smooth vector field, we may associate a $1$-form $\omega^1_{\bm{X}}$ with $\bm{X}$ via the Riemannian metric $g$ by setting $\omega^1_{\bm{X}}:=g(\bm{X},\cdot)$, which gives rise to an isomorphism between the spaces of smooth vector fields and (smooth) $1$-forms. The divergence of $\bm{X}$, denoted $\text{div}(\bm{X})$, is given by $\text{div}(\bm{X}):=\star d\star \omega^1_{\bm{X}}$, where $\star$ denotes the Hodge star operator and $d$ the exterior derivative, while the curl of $\bm{X}$, denoted $\text{curl}(\bm{X})$, is the unique vector field satisfying $\omega^1_{\text{curl}(\bm{X})}=\star d\omega^1_{\bm{X}}$. We call $\bm{X}\in \mathcal{V}^{\omega}(M)$ a \textit{Beltrami field} if $\text{div}(\bm{X})=0$ and if there exists a smooth function $\lambda:M\rightarrow \mathbb{R}$ with $\text{curl}(\bm{X})=\lambda \bm{X}$. Note that we only require smoothness of the metric and proportionality function, while we assume $\bm{X}$ to be real analytic. If $g$ and $\lambda$ are real analytic, then every smooth vector field satisfying $\text{curl}(\bm{X})=\lambda \bm{X}$ and $\text{div}(\bm{X})=0$ is, by standard elliptic estimates, necessarily real analytic. We use the words zero set and nodal set synonymously throughout the text.
\newline
\newline
Before stating our main theorem let us give two definitions used therein
\begin{defn}[Real analytic $1$-(sub)manifold with $C^1$-endpoints]
\label{MD1}
Suppose $M$ is a $3$-manifold. Let $L\subset M$ be a (not necessarily connected) subset of $M$ which is equipped with a $C^1$-atlas, turning it into a $1$-dimensional manifold with (possibly empty) boundary, such that the transition functions are all $C^1$-diffeomorphisms and such that the restrictions of the transition functions to the manifold interior of $L$ are real analytic diffeomorphisms. We say that $L$ is a \textit{real analytically embedded $1$-submanifold with $C^1$-endpoints}, if $L$ is a $C^1$-embedded $1$-submanifold of $M$ and if the manifold interior $\text{int}(L)$, with respect to the induced real analytic structure, is a real analytically embedded submanifold of $M$. We say the endpoints are non-empty if the manifold boundary of $L$ is non-empty and we say the endpoints are non-present otherwise.
\end{defn}
\begin{defn}[Tame knots]
\label{MD2}
Let $M$ be a $3$-manifold. A subset $K\subset M$ is called a \textit{knot} if there exists a homeomorphism $f:S^1\rightarrow K$. Now let $M$ be equipped with a smooth Riemannian metric $g$, then we call a knot $K\subset M$ \textit{tame} if there exists a continuous map $\Gamma:[0,l]\rightarrow M$ for some $l>0$ satisfying the following:
\begin{itemize}
\item $\Gamma([0,l])=K$.
\item $\Gamma(0)=\Gamma(l)$ and $\Gamma|_{[0,l)}$ is injective.
\item $\Gamma\in C^1([0,l],M)$ and $\Gamma|_{(0,l)}\in C^{\infty}((0,l),M)$.
\item $|\dot{\Gamma}(s)|_g=1$ for all $s\in [0,l]$.
\item $0\leq \int_0^l\kappa_{\Gamma}(s)ds<+\infty$, where $\kappa_{\Gamma}$ denotes the geodesic curvature of $\Gamma$.
\end{itemize}
We refer to $\int_0^l\kappa_{\Gamma}(s)ds$ as the \textit{total geodesic curvature}.
\end{defn}
\textit{Remark:} For the special case $M=\Omega\subseteq \mathbb{R}^3$ being an open subset of standard Euclidean space, equipped with the Euclidean metric, one easily checks that any tame knot in the sense of \cref{MD2} is of finite length and its unit tangent is of bounded variation. Hence if such a tame knot according to our definition is viewed as a subset of $\mathbb{R}^3$, it follows from \cite{M50}, \cite{Sul08} that it is also tame in the classical sense, i.e. there exists an ambient isotopy of Euclidean $3$-space transforming $K$ to a polygonal knot, or equivalently into a smoothly embedded circle. The idea to use bounded total curvature to obtain tame knots, in the Euclidean setting, was introduced in \cite{M50}. Here we adapted this notion to our setting.
\begin{thm}[Main theorem, Structure of nodal sets of real analytic Beltrami fields]
\label{MT3}
Let $M$ be a $3$-manifold which is equipped with a smooth Riemannian metric. Suppose $\bm{X}\in \mathcal{V}^{\omega}(M)$ is a Beltrami field, which is not the zero vector field and define $K:=\{p\in M|\bm{X}(p)=0 \}$. Then the Hausdorff dimension of $K$ is either $0$ or $1$ and there exists a locally finite, countable family of disjoint sets $\{A,L_1,L_2,\dots\}$ such that
\begin{itemize}
\item $K=A\sqcup \bigsqcup_{n\in \mathbb{N}}L_n$, where $\sqcup$ indicates that the union is disjoint.
\item The set $A$ is closed in $M$ and is either empty or consists of isolated points.
\item Each $L_n$ is either empty or a (non-empty) analytically embedded, connected, $1$-submanifold without boundary. Further each $L_n$ satisfies exactly one of the following conditions
\begin{itemize}
\item $\text{clos}(L_n)=L_n$, i.e. $L_n$ is closed in $M$.
\item $\text{clos}(L_n)$ is a connected, real analytically embedded $1$-submanifold with non-empty $C^1$-endpoints.
\item $\text{clos}(L_n)\setminus L_n$ is non-empty and $\text{clos}(L_n)$ is a tame knot.
\end{itemize}
\end{itemize}
If we let further $\mathcal{I}\subseteq A$ denote the isolated points of $K$, then $K\setminus \mathcal{I}=\bigcup_{n\in \mathbb{N}}\text{clos}(L_n)$.
\end{thm}
Let us state some simple implications
\begin{cor}
\label{MC4}
Let $M$ be a compact $3$-manifold which is equipped with some smooth metric. Suppose $\bm{X}\in \mathcal{V}^{\omega}(M)$ is a Beltrami field, which is not the zero vector field and let $K$ denote the zero set of $\bm{X}$ and $\mathcal{I}$ denote the isolated points of $K$. Then $\mathcal{I}$ consists of at most finitely many points and exactly one of the following two situations occurs
\begin{itemize}
\item $K\setminus \mathcal{I}=\emptyset$.
\item $K\setminus \mathcal{I}$ is non-empty and a finite union of tame knots and real analytically embedded $1$-submanifolds with $C^1$-endpoints ($C^1$-)diffeomorphic to $[0,1]$ which intersect in at most finitely many points.
\end{itemize}
\end{cor}
\begin{cor}
\label{MC5}
Let $\Omega\subset \mathbb{R}^3$ be a bounded domain of the standard Euclidean $3$-space (equipped with the standard metric) and let $\bm{X}$ be a real analytic Beltrami field on $\Omega$, which is not the zero vector field. Suppose that $\text{dist}\left(\partial \Omega,K \right)>0$, where $K$ is the zero set of $\bm{X}$ and the distance is the usual Euclidean distance. Then after removing at most finitely many isolated points, the set $K$ is either empty or is the finite union of tame knots and real analytically embedded $1$-submanifolds with $C^1$-endpoints, diffeomorphic to $[0,1]$, which intersect in at most finitely many points.
\end{cor}
As mentioned before if we view the tame knots in \cref{MC5} as knots in $\mathbb{R}^3$, then they are also tame in the classical sense. Observe also that we do not make any kind of regularity assumptions on the boundary of $\Omega$. A more intuitive way of phrasing the results of \cref{MC4} and \cref{MC5} is to say that the zero sets consist of at most finitely many isolated points and a finite collection of closed ($\cong S^1$) and open ($\cong [0,1]$) (well-behaved) strings, which intersect each other in at most finitely many points.
\newline
Now let us turn to the second part of our results. To this end we make the following definition
\begin{defn}[Torus knot, (Open) solid toroidal annulus]
\label{MD6}
Let $(p,q)\in \mathbb{N}^2$ be (strictly) positive integers which are coprime, then we define the map
\begin{equation}
\label{M1}
T_{p,q}:\mathbb{R}\rightarrow \mathbb{R}^3, t\mapsto \left(\cos(qt)(2+\cos(pt)),\sin(qt)(2+\cos(pt)),\sin(pt) \right)
\end{equation}
and call the image $\mathcal{T}_{p,q}:=T_{p,q}(\mathbb{R})$ the $(p,q)$-\textit{torus knot}.
\newline
Further we define
\begin{equation}
\label{M2}
\mathcal{T}_A:=\left\{(x,y,z)\in \mathbb{R}^3\left|\frac{1}{4}<\left(\sqrt{x^2+y^2}-2 \right)^2+z^2<\frac{9}{4} \right\}\right.,
\end{equation}
which we call \textit{(open) solid toroidal annulus}.
\end{defn}
See \cref{Figure1} in \cref{PP7E1} at the end of the paper for an illustration of $\mathcal{T}_{2,3}$ and $\mathcal{T}_A$.
\newline
Our main result regarding the possible complexity of nodal sets of real analytic Beltrami fields is the following.
\begin{thm}
\label{MP7}
Let $\mathcal{T}_A$ be the solid toroidal annulus, then given any (strictly) positive, coprime integers $(p,q)\in \mathbb{N}^2$ there exist countable infinitely many distinct real analytic metrics $(g_{p,q,k})_{k\in \mathbb{Z}}$ on $\mathcal{T}_A$ and countable infinitely many distinct, real analytic vector fields $(\bm{X}_{p,q,k})_{k\in \mathbb{Z}}$ on $\mathcal{T}_A$ such that
\begin{equation}
\label{M3}
\text{curl}_{g_{p,q,k}}\left(\bm{X}_{p,q,k}\right)=\bm{X}_{p,q,k}\text{ for all }k\in \mathbb{Z}
\end{equation}
with respect to the standard orientation on $\mathcal{T}_A$ and such that
\begin{equation}
\label{M4}
\{x\in \mathcal{T}_A| \bm{X}_{p,q,k}(x)=0 \}=\mathcal{T}_{p,q}\text{ for all }k\in \mathbb{Z}.
\end{equation}
The metrics $g_{p,q,k}$ and vector fields $\bm{X}_{p,q,k}$ can be explicitly expressed in terms of Cartesian coordinates by means of elementary functions alone.
\end{thm}
\section{Proof of \cref{MT3}}
First let us recall that the order of a zero $p$ of a given smooth vector field $\bm{X}$ is defined, after fixing any chart $\mu$ around $p$, as the minimum of the orders of the zero $\mu(p)$ of the corresponding local expressions $X^j\circ \mu^{-1}$ of $\bm{X}$. This definition is independent of the choice of chart. We denote the order of a given zero $p$ by $\Omega(p)$.
\newline
The proof consists of several steps, which we formulate as lemmas to increase readability
\begin{lem}
\label{Step1}
Suppose we are in the setting of \cref{MT3}, then there exists a countable, disjoint family of sets $\{A,L_1,L_2,\dots\}$ such that $K=A\sqcup\bigsqcup_{n\in \mathbb{N}}L_n$ and such that each of the $L_n$ is either empty or an analytically embedded $1$-submanifold without boundary.
\end{lem}
\underline{Proof of \cref{Step1}:} Given $n\in \mathbb{N}$ let $S_n:=\{p\in K| \Omega(p)=n \}$ and observe that by analyticity of $\bm{X}$ and since we assume it not to be the zero vector field
\begin{equation}
\label{PT31}
K=\bigsqcup_{n\in \mathbb{N}}S_n,
\end{equation}
where $\bigsqcup$ indicates that the union is disjoint. Fix any non-empty $S_n$. We claim that either $S_n$ consists entirely of isolated points or is a union of isolated points and a (non-empty and not necessarily connected) real analytically embedded $1$-submanifold without boundary. Fix some $p\in S_n$ and a coordinate chart $\mu_p:U_p\rightarrow V_p\subseteq \mathbb{R}^3$ around $p$ with $\mu(p)=0$ (one may choose normal coordinates to simplify some calculations). By definition of $S_n$ there exists a multi-index $|\beta|=n-1$ and some $X^i$ with $\left((\partial_1\partial^{\beta}X^i)(p),(\partial_2\partial^{\beta}X^i)(p),(\partial_3\partial^{\beta}X^i)(p) \right)\neq 0$ and after possibly shrinking $U_p$ the gradient of this function never vanishes within $U_p$. For simplicity assume $i=1$. We can then define the vector field
\begin{equation}
\label{PT32}
\bm{h}:=\left((\partial^{\beta}X^j)\circ \mu^{-1}_p\right)e_j:V_p\rightarrow \mathbb{R}^3,
\end{equation}
where $e_j$ denote the standard basis vectors of $\mathbb{R}^3$. Consider the Jacobian
\begin{equation}
\label{PT33}
M_p:=(D\bm{h})(0).
\end{equation}
By choice of $\bm{h}$ we certainly have $\text{rank}(M_p)\in \{1,2,3\}$. On the other hand by applying $\partial^{\beta}$ to both sides of the local expressions of $\text{curl}(\bm{X})=\lambda \bm{X}$ and $\text{div}(\bm{X})=0$, keeping in mind the definition of $S_n$, one obtains
\begin{equation}
\label{PT34}
(\partial_ih^j)(0)=(\partial_jh^i)(0)\text{ for all }1\leq i,j\leq 3\text{ and }(\partial_ih^i)(0)=0,
\end{equation}
where in the latter equation we use Einstein's summation convention. The relations in (\ref{PT34}) contradict the assumption $\text{rank}(M_p)=1$, so that we must have $\text{rank}(M_p)\in \{2,3\}$. The details of this argument are carried out in a paper, whose preprint version can be found on arXiv\footnote{arXiv identifier: 2005.06590} [p. 11, proof of 2nd part of proposition 2.8]. If the rank is $3$, then the inverse function theorem implies that, after possibly shrinking $U_p$, $0\in V_p$ is the unique solution of $\bm{h}(x)=0$ with $x\in V_p$. It follows from the definition of $S_n$ that $S_n\cap U_p=\{p\}$, i.e. $p$ is an isolated point in this case. Now $S_n$ may consist of isolated points alone, then we are done. Thus let from now on $\mathcal{L}_n\neq \emptyset$ denote the non-isolated points of $S_n$. Fix any $p\in \mathcal{L}_n$, then by our previous arguments we must have $\text{rank}(M_p)=2$. We may assume that $(\nabla h^1)(0)$ and $(\nabla h^2)(0)$ are linearly independent, where $\nabla$ denotes the Euclidean gradient. Define $\hat{\bm{h}}:V_p\rightarrow \mathbb{R}^2$ $x\mapsto (h^1(x),h^2(x))$, then $D\hat{\bm{h}}(0)$ has rank two and thus, after possibly interchanging the role of the coordinate axis, $(\partial_2\hat{\bm{h}})(0)$ and $(\partial_3\hat{\bm{h}})(0)$ are linearly independent. It then follows from the real analytic implicit function theorem \cite[Theorem 2.3.5]{KrPa02} that there exists an open interval $0\in I\subseteq \mathbb{R}$ and open subset $0\in W\subseteq \mathbb{R}^2$ with $I\times W\subseteq V_p$ and a real analytic function $\phi=(\phi_1,\phi_2):I\rightarrow W$ with $\phi(0)=0$ and
\begin{equation}
\label{PT35}
\forall\text{ }(t,x)\in I\times W\text{: }\hat{\bm{h}}(t,x)=0\Leftrightarrow x=\phi(t).
\end{equation}
After shrinking $U_p$ if necessary, we obtain $V_p=I\times W$. Further, since $p\in S_n$ is not isolated, there exists a sequence $(p_k)_k\subseteq S_n\setminus \{p\}$ converging to $p$. Thus for all high enough indices we have $p_k\in S_n\cap U_p$ and by definition of $S_n$ and due to (\ref{PT35}) $\mu_p(p_k)=(t_k,\phi(t_k))$. Observe that $t_k\neq 0$ for all such $k$ since otherwise $\mu_p(p_k)=(0,\phi(0))=(0,0)=\mu_p(p)$ and hence $p_k=p$, a contradiction. But since $(p_k)_k$ converges to $p$ and $\mu_p(p)=0$, there is a sequence $(t_k)_k\subset I\setminus \{0\}$ converging to $0$ with $\mu_p(p_k)=(t_k,\phi(t_k))$. Now fix any multi-index $|\alpha|\leq n-1$ and define the functions
\begin{equation}
\label{PT36}
f^j_{\alpha}:=(\partial^\alpha X^j)\circ \mu^{-1}_p\circ \left(Id|_I\times \phi \right):I\rightarrow \mathbb{R}.
\end{equation}
Observe that the $f^j_{\alpha}$ are real analytic as compositions of real analytic functions. Further we have by definition of $S_n$, $f^j_{\alpha}(0)=(\partial^\alpha X^j)(p)=0$ and $f^j_{\alpha}(t_k)=(\partial^\alpha X^j)(p_k)=0$ for all $k$. We conclude that the set $\{t\in I| f^j_{\alpha}(t)=0\}$ has an accumulation point and since the $f^j_{\alpha}$ are real analytic and $I$ is an interval, we find $f^j_{\alpha}\equiv 0$ on $I$ for every $1\leq j \leq 3 $ and multi-index $|\alpha|\leq n-1$. Now consider the set
\begin{equation}
\label{PT37}
l_p:=\{q\in U_p| (\partial^{\beta}X^1)(q)=0=(\partial^{\beta}X^2)(q) \}\end{equation}
and observe that $q\in l_p$ if and only if $\hat{\bm{h}}(\mu_p(q))=0$. We then conclude from (\ref{PT35}) that $\mu_p(q)=(t,\phi(t))$ for some suitable $t\in I$. Let $|\alpha|\leq n-1$ be any multi-index, then
\[
(\partial^{\alpha}X^j)(q)=((\partial^{\alpha}X^j)\circ{\mu^{-1}_p})(\mu_p(q))=((\partial^{\alpha}X^j)\circ{\mu^{-1}_p})(t,\phi(t))=f^j_{\alpha}(t)=0
\]
by our previous findings. Hence $\Omega(q)\geq n$ and by choice of $U_p$ we know that $\nabla\left(\partial^{\beta}X^1 \right)(q)\neq 0$ for every $q\in U_p$. Thus we must have $\Omega(q)=n$ and overall we obtain $l_p\subseteq S_n\cap U_p$. The converse implication follows trivially from the definition of $S_n$ and $l_p$, so that we arrive at
\begin{equation}
\label{PT38}
S_n\cap U_p=l_p\text{ for some open neighbourhood }U_p\text{ around }p.
\end{equation}
We can finally define $\psi_p:I\rightarrow l_p, t\mapsto \mu^{-1}_p((t,\phi(t)))$, which gives rise to a homeomorphism between $I$ and the open subset $l_p=S_n\cap U_p$ of $S_n$, where the inverse is given by $\pi_1\circ \mu_p$, with $\pi_1$ being the projection onto the first component. Observe that since $\psi_p$ is a homeomorphism, none of the points in $U_p\cap S_n$ are isolated points and so in fact we have $U_p\cap S_n=U_p\cap \mathcal{L}_n$, so that the maps $\psi_p$ give rise to an atlas of $\mathcal{L}_n$. It is easy to check that the transition functions are real analytic and that the so obtained real analytic manifold $\mathcal{L}_n$ is in fact real analytically embedded.
\newline
To conclude the proof of the first lemma, let $A_n$ denote the set of isolated points of $S_n$, then $S_n=A_n\sqcup \mathcal{L}_n$. We can further decompose $\mathcal{L}_n$ into its connected components $\mathcal{L}_n=\bigsqcup_{m\in \mathbb{N}}C_{m,n}$. Define $A:=\bigsqcup_{n\in \mathbb{N}}A_n$ and identify the $(L_n)$'s with connected components $C_{n,m}$ of the $\mathcal{L}_n$, then this gives us our countable, disjoint family with the claimed properties. $\square$
\begin{lem}
\label{Step2}
The family of sets constructed in \cref{Step1} is locally finite and the set $A$ consists of isolated points alone and is a closed subset of $M$.
\end{lem}
\underline{Proof of \cref{Step2}:} We will first prove that this family is locally finite. First note that $K\subset M$ is closed and therefore we may choose for any $p\in M\setminus K$ the set $M\setminus K$ as an open neighbourhood not intersecting any family member. Thus let $p\in K$. It follows immediately from definition that the order of a zero is locally nonincreasing, i.e. we can find an open neighbourhood $U_p$ around $p$ such that $\Omega(q)\leq \Omega(p)$ for all $q\in U_p\cap K$. In other words $U_p\cap K\cap S_n=\emptyset$ for every $n>\Omega(p)$. Now fix any $1\leq n\leq \Omega(p)$ and observe that by definition of the sets $S_n$, they are all semianalytic (see \cite{BM88} for a thorough introduction into this topic). It then follows from \cite[Corollary 2.7]{BM88} that the family of connected components of the $S_n$ are locally finite. But $A_n$ consists exactly of the isolated points of $S_n$, while we had shown in the previous proof that for each $q\in \mathcal{L}_n$ there is an open set $V_q$ such that $V_q\cap S_n=V_q\cap \mathcal{L}_n$. This implies that there exists an open set $V$ of $M$ with $S_n\cap V= \mathcal{L}_n$ and therefore the connected components of $\mathcal{L}_n$ give rise to connected components of $S_n$. Hence the family of connected components of $\mathcal{L}_n$ is locally finite, i.e. we can find an open neighbourhood $U_n$ around $p$ such that at most finitely many of the connected components of $\mathcal{L}_n$ intersect $U_n$. Define $U:=U_p\cap \bigcap_{n=1}^{\Omega(p)}U_n$, then $U$ is an open neighbourhood around $p$ which intersects only finitely many of the connected components of the $\mathcal{L}_n$ for $1\leq n\leq \Omega(p)$ and hence also only at most finitely many of the $L_m$. Since $p\in K$ was arbitrary the family is indeed locally finite.
\newline
\newline
Now consider the set $A$ and suppose that $(p_k)_k\subseteq A$ is a sequence converging to some $p\in M$. We claim that there exists a subsequence of $(p_k)_k$ which is entirely contained in some $A_n$ and that $p\in A_n$ for the same $n$. To this end we observe that $p\in K$ by closedness of $K$. Just like before we can find an open neighbourhood $U_p$ of $p$ such that $A\cap U_p=\bigsqcup_{n=1}^{\Omega(p)}A_n\cap U_p$ and thus at least one of these finitely many $A_n$ must contain infinitely many members of the sequence $(p_k)_k$. Fix any $n$ with this property and denote the corresponding subsequence again by $(p_k)_k$. Now let  $\mathring{S}^{(d)}_n$ for $0\leq d\leq 3$ denote the subsets of $S_n$ which admit an open neighbourhood $U$ in $M$ such that $U\cap S_n$ is an analytically embedded $d$-dimensional submanifold of $M$, then by \cite[Theorem 7.2, Remark 7.3]{BM88} and \cite[Theorem 6.11]{MaTr07} all these sets are semianalytic. By our findings we have $\mathring{S}^{(3)}_n=\emptyset$, $\mathring{S}^{(2)}_n=\emptyset$, $\mathring{S}^{(1)}_n=\mathcal{L}_n$ and $\mathring{S}^{(0)}_n=A_n$. Therefore $A_n$ is semianalytic and $p\in \text{clos}(A_n)$. Hence by the curve selection lemma, \cite[Lemma 2.1]{Loo84}, \cite[Lemma 6.6]{MaTr07}, we can in particular find a continuous curve $\gamma:[0,\delta)\rightarrow M$ for some $\delta>0$ with $\gamma(0)=p$ and $\gamma((0,\delta))\subseteq A_n$. Since $\gamma((0,\delta))$ is connected and $A_n$ consists precisely of the  isolated points of $S_n$ the image $\gamma((0,\delta))$ is a single point. By continuity of $\gamma$ this point must coincide with $p$ and hence $p\in A_n$, which proves the claim. In particular $p\in A$ and hence $A$ is a closed subset of $M$. If we assume that $p\in A$ is not an isolated point in $A$ we can find a sequence in $A$ converging to $p$ and consisting of distinct elements. But as we have seen we can extract a subsequence and find some $n$ such that the limit point $p$ and the subsequence are contained in $A_n$. But since $A_n$ consists only of isolated points by definition all but finitely many elements of the sequence must be equal to $p$, a contradiction. This concludes the proof. $\square$
\begin{lem}
\label{Step3}
Suppose we are in the setting of \cref{MT3} and let $\mathcal{I}$ be the set of all isolated points of $K$, then $K\setminus \mathcal{I}=\bigcup_{n\in \mathbb{N}}\text{clos}(L_n)$, where the $L_n$ are the sets constructed in \cref{Step1}.
\end{lem} 
\underline{Proof of \cref{Step3}:}
Every $L_n$ is an analytically embedded $1$-submanifold and so no $L_n$ contains an isolated point, i.e. $L_n\subseteq K\setminus \mathcal{I}$ and since $K$ is closed we also have $\text{clos}(L_n)\subseteq K$. We then in addition have $\text{clos}(L_n)\setminus L_n\subseteq K\setminus\mathcal{I}$ and therefore $\text{clos}(L_n)\subseteq K\setminus\mathcal{I}$ for every $n$. For the converse implication fix any $p\in K\setminus\mathcal{I}$. By \cref{Step2} we know that the family $\{A,L_1,L_2\dots\}$ is locally finite. Hence fix some open neighbourhood $U_p$ around $p$ which only intersects finitely many of the $L_n$. Since $p$ is not an isolated point of $K$ we can find a sequence $(p_k)_k$, contained in $(K\cap U_p)\setminus \{p\}$, which converges to $p$. Now if this sequence were to contain a subsequence contained in $A$, then by closedness of $A$ we would have $p\in A$ and since $A$ consists only of isolated points that would imply that the corresponding subsequence is a constant sequence, except for at most finitely many elements. This contradicts the choice of the original sequence. Therefore, after removing at most finitely many elements if necessary, the sequence $(p_k)_k$ is contained in the union of finitely many $L_n$. Thus there must exist at least one $L_n$ containing a subsequence of $(p_k)_k$ which converges to $p$, i.e. $p\in \text{clos}(L_n)\subseteq \bigcup_{n\in \mathbb{N}}\text{clos}(L_n)$. $\square$
\begin{lem}
\label{Step4}
The sets $L_n$ constructed in \cref{Step1} satisfy the following: If $\text{clos}(L_n)\setminus L_n$ is non-empty, then it contains either exactly $1$ or exactly $2$ elements.
\end{lem} 
\underline{Proof of \cref{Step4}:} Fix any $L_n$ with $\text{clos}(L_n)\setminus L_n\neq \emptyset$. Since $M$ is Hausdorff and since compact subsets of Hausdorff spaces are closed, we know that $L_n$ is in particular not compact. On the other hand we know that $L_n$ is an analytically embedded $1$-manifold without boundary. Thus by the classification of $1$-manifolds \cite{M65} we see that $L_n$ is diffeomorphic to $(0,1)$. Let $\psi:(0,1)\rightarrow L_n$, denote any fixed diffeomorphism. Now fix any $p\in \text{clos}(L_n)\setminus L_n$ and let $(p_k)_k\subset L_n$ be any sequence in $L_n$ converging to $p$. By definition of $\psi$ there is a sequence $(t_k)_k$ in the open unit interval with $p_k=\psi(t_k)$. After choosing a subsequence of $(t_k)_k$, we may assume that it converges to some $t$ within $[0,1]$ and since $p$ is not contained in $L_n$, we must have $t\in \{0,1\}$. We assume from now $t=1$, since the other case can be treated identically. We claim that if $(s_k)_k\subseteq (0,1)$ is any other sequence converging to $1$, then $\psi(s_k)\subset M$ converges to $p$ as well. Once this is shown \cref{Step4} will be proven. So let $(s_k)_k$ be any such fixed sequence and let $(s_{k_m})_m$ be any fixed subsequence of $(s_k)_k$. We will show that each such subsequence admits yet another subsequence such that $\psi(s_{k_{m_j}})_j$ converges to $p$. Then the claim will follow. To simplify notation we write $(s_k)_k$ instead of $(s_{k_m})_m$. We already have argued before that the $\mathcal{L}_n$ are semianalytic and hence all of their connected components are semianalytic \cite[Corollary 2.7]{BM88}. From this one concludes that $\psi\left((0,\frac{1}{2})\right)$ and $\psi\left((\frac{1}{2},1)\right)$ are both semianalytic subsets of $M$. Now $p\in \text{clos}\left( \psi\left((\frac{1}{2},1)\right)\right)$ and so by the curve selection lemma we can find some $\delta>0$ and a real analytic curve $\gamma:(-\delta,\delta)\rightarrow M$ with $\gamma(0)=p$ and $\gamma((0,\delta))\subseteq \psi\left((\frac{1}{2},1)\right)$. This allows us to define the function
\begin{equation}
\label{PT39}
\phi:\left(0,\frac{\delta}{2}\right)\rightarrow \left(\frac{1}{2},1\right),t\mapsto \psi^{-1}(\gamma(t)).
\end{equation}
We claim that for any sequence $(\tau_l)_l\subset \left(0,\frac{\delta}{2}\right)$, converging to $0$, we can extract a subsequence $(\tau_{l_m})_m$ such that $\phi(\tau_{l_m})$ converges to $1$. So let $(\tau_l)_l$ be such a sequence, then $\phi(\tau_l)$ is bounded and hence we can extract a subsequence (again denoted $(\tau_l)_l$) such that $\phi(\tau_l)$ converges to some element in $\left[\frac{1}{2},1\right]$. But we have the equality $\gamma(\tau_l)=\psi(\phi(\tau_l))$. By properties of $\gamma$ the former converges to $p$, which is not an element of $L_n$ and thus $\phi(\tau_l)$ must converge to $1$. This proves the claim. From this we in particular conclude that if we let $J:=\phi\left(\left(0,\frac{\delta}{2}\right)\right)\subseteq \left(\frac{1}{2},1\right)$, then by continuity $J$ must be an interval and by what we have shown $1\in \text{clos}(J)$, i.e. $J$ is a non-degenerate interval. Now fix any null-sequence $(\tau_l)_l$ for which $\phi(\tau_l)$ converges to $1$. We recall that $(s_k)_k$ was converging to $1$ and so does any of its subsequences. So we can find for $\phi(\tau_1)$ some $k_1$ with $\phi(\tau_1)<s_{k_1}$ since both sequences converge to $1$ from below. But then we can also find some $\tau_{l_2}$ with $s_{k_1}<\phi(\tau_{l_2})$. That way we may select subsequences of $(\tau_l)_l$ and $(s_k)_k$, again denoted in the same way, with $\phi(\tau_1)<s_1<\phi(\tau_2)<s_2<\dots$. Now observe that $\phi(\tau_k)\in J$ for every $k$ and that $J$ is a non-degenerate interval, i.e. $[\phi(\tau_k),\phi(\tau_{k+1})]\subseteq J$ for every $k$. By construction of our subsequence we have $s_k\in  [\phi(\tau_k),\phi(\tau_{k+1})]\subseteq J$, i.e. $s_k\in J$ for every $k$, and by definition of $J$ there exists for every $k$ a $\sigma_k\in \left(0,\frac{\delta}{2} \right)$ with $s_k=\phi(\sigma_k)$. We claim that after extracting a subsequence, if necessary, the $\sigma_k$ converge to $0$. Since $(\sigma_k)_k\subset \left(0,\frac{\delta}{2}\right)$ we can select any convergent subsequence (again denoted $\sigma_k$) and observe that $\gamma(\sigma_k)=\psi(\phi(\sigma_k))=\psi(s_k)$ $\Leftrightarrow s_k=\psi^{-1}\left(\gamma(\sigma_k)\right)$. Since $\gamma\left(\left(0,\delta \right)\right)\subset L_n$, we see that if $\sigma_k$ converges to some $\sigma\in \left(\left.0,\frac{\delta}{2}\right.\right]$, then $(s_k)_k$ must converge to some element in $(0,1)$, a contradiction. Thus $(\sigma_k)_k$ converges to $0$ and in conclusion $\psi(s_k)=\gamma(\sigma_k)\rightarrow \gamma(0)=p$ by properties of $\gamma$. $\square$
\begin{lem}
\label{Step5}
For each of the sets $L_n$, constructed in \cref{Step1}, the closure $\text{clos}(L_n)$ is either a real analytically embedded $1$-manifold with $C^1$-endpoints or a knot.
\end{lem}
\underline{Proof of \cref{Step5}:} If $L_n$ is closed in $M$, the statement follows from \cref{Step1}. Thus, according to \cref{Step4}, let $\text{clos}(L_n)\setminus L_n=\{p_0,p_1\}$, where we allow $p_0=p_1$. As argued in the previous proof there is a diffeomorphism $\psi:(0,1)\rightarrow L_n$. If $p_0\neq p_1$ we can label them such that there are sequences $(t_k)_k$, $(s_k)_k$ in $(0,1)$ converging to $0$ and $1$ respectively with $\psi(t_k)$ and $\psi(s_k)$ converging to $p_0$ and $p_1$ respectively. In case of $p_0=p_1\equiv p$, we may (after replacing $\psi$ by $\psi(1-\cdot)$ if necessary) assume that there exists a sequence $(s_k)_k$ in $(0,1)$ converging to $1$ with $\psi(s_k)$ converging to $p$. Now define $C_0:=\psi\left(\left(0,\frac{1}{2}\right)\right)$ and $C_1:=\psi\left(\left(\frac{1}{2},1\right)\right)$, which are semianalytic sets, and observe that $p_0$ and $p_1$ are contained in the respective closures if the points are distinct or are contained in the closure of $C_1$ if they coincide. Hence by the curve selection lemma we may define the following real analytic curves
\begin{equation}
\label{PT310}
\gamma_{0 \slash 1}: \left(-\delta,\delta\right)\rightarrow M, \text{ }\gamma_{0 \slash 1}(0)=p_{0 \slash 1}\text{ and }\gamma_{0 \slash 1}((0,\delta))\subseteq C_{0 \slash 1}.
\end{equation}
Observe that after possibly shrinking $\delta$ we may assume that $\dot{\gamma}_{0\slash 1}(t)\neq 0$ for all $0<|t|<\delta$, since otherwise the real analyticity of the $\gamma_{0\slash 1}$ implies that the curves must be constant, contradicting the fact that $p_{0\slash 1}$ is not contained in $L_n$. In a similar spirit one can argue that there must be some $\delta>0$ such that $\gamma_{0 \slash 1}|_{[0,\delta)}$ is injective. We observe that this implies that the maps
\begin{equation}
\label{PT311}
\gamma_{0 \slash 1}:\left[0,\frac{\delta}{2} \right]\rightarrow \gamma_{0 \slash 1}\left(\left[0,\frac{\delta}{2} \right]\right),
\end{equation}
are continuous bijections from a compact space into a Hausdorff space and therefore are homeomorphisms. In particular the restrictions $\gamma_{0 \slash 1}|_{\left[0,\frac{\delta}{2}\right)}$ are homeomorphisms onto their images.
\newline
\newline
\textit{Case 1:} Assume there does not exist a sequence $(t_k)_k$ in $(0,1)$ converging to $0$, such that $\psi(t_k)$ converges to $p_1$. We claim that under this assumption, the set $\gamma_1\left(\left[0,\frac{\delta}{2}\right)\right)$ is an open subset of $\text{clos}(L_n)$. Observe that if $\text{clos}(L_n)\setminus L_n$ consists of two points, then the assumption is necessarily satisfied. Further note that the images of the maps in (\ref{PT311}) are disjoint since $C_0$ and $C_1$ are disjoint. So for the upcoming arguments it is enough to consider $\gamma\equiv \gamma_1$, since identical arguments apply to $\gamma_0$.
\newline
We observe first that $\gamma:(0,\delta)\rightarrow L_n$ is smooth and by choice of $\delta$ its derivative is nowhere vanishing. Thus the inverse function theorem implies that $\gamma$ is an open map and in particular $\gamma\left(\left(0,\frac{\delta}{2}\right)\right)$ is an open subset of $L_n$. In addition $\text{clos}(L_n)\setminus L_n$ consists of finitely many points and hence is a closed subset of $\text{clos}(L_n)$, i.e. $L_n$ is an open subset of $\text{clos}(L_n)$ and thus overall $\gamma\left(\left(0,\frac{\delta}{2}\right)\right)$ is an open subset of $\text{clos}(L_n)$. On the other hand consider the following map
\[
\Psi:(0,1]\rightarrow \text{clos}(L_n), t\mapsto \begin{cases}
\psi(t) & 0<t<1 \\
p_1 &t=1
\end{cases}
\]
and observe that for any $0<\epsilon<1$ we have the equality
\begin{equation}
\label{PT312}
\Psi((\epsilon,1])=\begin{cases}
\text{clos}(L_n)\setminus\left(\psi((0,\epsilon])\cup \{p_0\} \right)& \text{ if }p_0\neq p_1 \\
\text{clos}(L_n)\setminus \psi((0,\epsilon])& \text{ }p_0=p_1
\end{cases}.
\end{equation}
It follows from similar arguments as in the proof of \cref{Step4} that $\left(\psi((0,\epsilon])\cup \{p_0\} \right)$ is closed in $\text{clos}(L_n)$ in the first case and, given our assumption, that $\psi((0,\epsilon])$ is closed in $\text{clos}(L_n)$ in the second case. In addition one can also argue similarly that for $0<\epsilon<1$ close enough to $1$ we have the inclusion $\Psi((\epsilon,1])\subseteq \gamma\left([0,\frac{\delta}{2}) \right)$. Now fix any such $\epsilon$ with this property and observe
\[
\gamma\left(\left[0,\frac{\delta}{2}\right)\right)=\Psi((\epsilon,1])\cup \gamma\left(\left(0,\frac{\delta}{2}\right)\right),
\]
where the latter is an open set of $\text{clos}(L_n)$ and where $\Psi((\epsilon,1])$ is open as the complement of a closed subset, (\ref{PT312}). Hence according to (\ref{PT311}) we see that $\gamma|_{\left[0,\frac{\delta}{2}\right)}$ defines a homeomorphism onto an open subset of $\text{clos}(L_n)$.
\newline
Recall that we fixed $\delta$ such that $\dot{\gamma}(t)\neq 0$ for all $0<t<\delta$, but did not exclude the possibility that $\dot{\gamma}(0)=0$. To address this, we can find by real analyticity of $\gamma$ a natural number $N\in \mathbb{N}$, such that the following map
\begin{equation}
\label{PT313}
\tilde{\gamma}:\left[0,\left(\frac{\delta}{2}\right)^N\right)\rightarrow M, s\mapsto \gamma(\sqrt[N]{s})
\end{equation}
is $C^1$ up to the boundary and satisfies $\dot{\tilde{\gamma}}(0)\neq 0$. By the chain rule we see that in fact $\dot{\tilde{\gamma}}$ never vanishes on $\left[0,(\frac{\delta}{2})^N\right)$ and that it still defines a homeomorphism onto an open subset. Note also that $\tilde{\gamma}$ is always real analytic away from zero. We can then fix $\tilde{\gamma}$ as a chart around $p_1$ and equip $L_n$ with the atlas constructed in the proof of \cref{Step1}. If $p_0\neq p_1$ we can find a chart around $p_0$ in exactly the same way as we did for $p_1$. Keeping in mind that $\dot{\tilde{\gamma}}(s)\neq 0$ for all $0\leq s< (\frac{\delta}{2})^N$ it is straightforward to confirm that $\text{clos}(L_n)$ equipped with these charts is a real analytically embedded $1$-manifold with $C^1$-endpoints whose manifold boundary coincides with $\text{clos}(L_n)\setminus L_n$.
\newline
\newline
\textit{Case 2:} Suppose there exists a sequence $(t_k)_k$ in $(0,1)$ converging to $0$ such that $\psi(t_k)$ converges to $p_1$. It then follows from the same arguments as in the proof of \cref{Step4} that we must have $p_0=p_1\equiv p$ and that the following map
\begin{equation}
\label{PT314}
\Psi:[0,1]\rightarrow M, t\mapsto \begin{cases}
\psi(t) & 0<t<1 \\
p & t=0,1
\end{cases}
\end{equation}
is continuous. It follows immediately from construction that $\Psi|_{[0,1)}$ is injective and hence it descends to a homeomorphism $f:S^1\rightarrow \Psi([0,1])=\text{clos}(L_n)$, which proves that $\text{clos}(L_n)$ in this case is indeed a knot. $\square$
\begin{lem}
\label{Step6}
In the situation of case 2 of the previous proof, the corresponding knot $\text{clos}(L_n)$ is tame in the sense of \cref{MD2}.
\end{lem} 
\underline{Proof of \cref{Step6}:} Observe that since we are in case 2 of the previous step, we have $\text{clos}(L_n)\setminus L_n=\{p\}$ and $p\in \text{clos}(C_0)\cap \text{clos}(C_1)$, see (\ref{PT310}). Thus we may define two real analytic curves $\gamma_{0 \slash 1}$ with the same properties as in (\ref{PT310}). We set $\gamma\equiv \gamma_1$ since the other case may be treated identically. By what we had shown, after choosing $\delta$ small enough, we know that $\gamma\left(\left(0,\frac{\delta}{2}\right)\right)$ is an open (and connected) subset of $L_n$. If we let again $\psi$ denote our fixed diffeomorphism from $(0,1)$ to $L_n$, we see that $\psi^{-1}\left(\gamma\left(\left(0,\frac{\delta}{2}\right)\right)\right)$ is a connected, non-empty and open subset of $\left(\frac{1}{2},1\right)$, since $\gamma$ maps into $C_1$, (\ref{PT310}), and hence is an open interval. If we just like in the proof of \cref{Step4} define the function $\phi$, \ref{PT39}, then $J:=\phi\left(\left(0,\frac{\delta}{2}\right)\right)$ is an open interval and by our findings in the proof of \cref{Step4} we know that $1\in \text{clos}(J)$. Therefore there exists some $\frac{1}{2}\leq \epsilon_1<1$ with $J=(\epsilon_1,1)$. Recall that $\delta$ was chosen so small that $\gamma|_{[0,\delta)}$ is injective. One can then argue in the same fashion as in the proof of \cref{Step4} that we must have $\gamma\left(\frac{\delta}{2}\right)=\psi(\epsilon_1)$. Thus after possibly shrinking $\delta$ a little bit further we may assume $\epsilon_1>\frac{1}{2}$. We can finally replace $\gamma$ by its regularised version $\tilde{\gamma}=\gamma(\sqrt[N_1]{\cdot})$ for a suitable $N_1\in \mathbb{N}$, where $\tilde{\gamma}$ is $C^1$- up to the boundary, injective, has nowhere vanishing differential and is real analytic on the interior $\left(0, \left(\frac{\delta}{2}\right)^{N_1}\right)$. An identical argument provides us with a corresponding map $\tilde{\gamma}_0$ with image $\tilde{\gamma}_0\left(\left(0,\left(\frac{\delta}{2}\right)^{N_0}\right)\right)=\psi((0,\epsilon_0)) $ for a suitable $0<\epsilon_0<\frac{1}{2}<\epsilon_1<1$, with $\tilde{\gamma}_0(0)=p$ and $\tilde{\gamma}_0\left(\left(\frac{\delta}{2}\right)^{N_0}\right)=\psi(\epsilon_0)$. Thus the following curve
\begin{equation}
\label{PT315}
\hat{\gamma}:[0,1]\rightarrow M, t\mapsto \begin{cases}
\tilde{\gamma}_0\left(\left(\frac{\delta}{2} \right)^{N_0}\frac{t}{\epsilon_0} \right) & 0\leq t \leq \epsilon_0 \\
\psi(t) & \epsilon_0\leq t \leq \epsilon_1\\
\tilde{\gamma}_1\left(\left(\frac{\delta}{2} \right)^{N_1}\frac{1-t}{1-\epsilon_1} \right) & \epsilon_1\leq t \leq 1
\end{cases}
\end{equation}
is continuous by the gluing lemma. By construction $\hat{\gamma}|_{[0,1)}$ is injective, $\hat{\gamma}(0)=p=\hat{\gamma}(1)$ and $\hat{\gamma}([0,1])=\text{clos}(L_n)$. Further it is piecewise $C^1$ and piecewise smooth away from $0$ and $1$. In order to prove that $\text{clos}(L_n)$ is tame, we will consider the arc-length parametrisation $\Gamma$ of $\hat{\gamma}$. Define the function $s:[0,1]\rightarrow [0,l], \tau\mapsto \int_0^{\tau}|\dot{\hat{\gamma}}(t)|_g dt$ with $l:=\int_0^{1}|\dot{\hat{\gamma}}(t)|_g dt<+\infty$. This function is strictly monotonically increasing and onto, hence has an inverse $s^{-1}:[0,l]\rightarrow [0,1]$ and we define $\Gamma:=\hat{\gamma}\circ s^{-1}:[0,l]\rightarrow M$. Observe that $\Gamma$ preserves the properties of $\hat{\gamma}$, i.e. $\Gamma|_{[0,l)}$ is injective, $\Gamma(0)=\Gamma(l)$, $\Gamma([0,l])=\text{clos}(L_n)$ and $\Gamma$ is piecewise $C^1$ and piecewise smooth away from $0$ and $l$. We define for convenience $s_0:=s(\epsilon_0)$ and $s_1:=s(\epsilon_1)$. By definition we have $|\dot{\Gamma}(\tau)|_g=1$ for all $0\leq \tau\leq l$, where this equality holds for both (a priori possibly distinct) limits at the points $s_0$ and $s_1$. We claim that $\Gamma|_{(0,l)}$ is smooth. It is obviously enough to establish smoothness in $s_0$ and $s_1$. We will show smoothness in $s_0$, since the other case can be treated identically. Define the following function
\[
\sigma:I:=(\epsilon_0-\epsilon,\epsilon_0+\epsilon)\rightarrow \mathbb{R}, t\mapsto \int_0^{\epsilon_0}|\dot{\hat{\gamma}}(\tau)|_gd\tau+\int_{\epsilon_0}^{t}|\dot{\psi}(\tau)|_gd\tau,
\]
where $\epsilon>0$ is so small that $I\subset (0,\epsilon_1)$. Observe that $\sigma$ is smooth and strictly increasing, hence has a smooth inverse $\sigma^{-1}:\tilde{I}\rightarrow I$, which is again strictly increasing and where $\tilde{I}=\sigma(I)$ is an open interval containing $\sigma(\epsilon_0)=s(\epsilon_0)=s_0$. By definition of $\sigma$ we therefore find $\sigma^{-1}(\tau)=s^{-1}(\tau)$ for all $\tau\geq s_0$. We define
\begin{equation}
\label{PT316}
\tilde{\Gamma}_+:=\psi\circ \sigma^{-1}:\tilde{I}\rightarrow M
\end{equation}
and observe that for all $\tau\in \tilde{I}$ with $\tau\geq s_0$ we have $\tilde{\Gamma}_+(\tau)=\psi(\sigma^{-1}(\tau))=(\psi\circ s^{-1})(\tau)=\Gamma(\tau)$, where we used the strict monotonicity, i.e. $s^{-1}(\tau)\geq s^{-1}(s_0)=\epsilon_0$. Thus if we write $\tilde{I}=(s_0-s_l,s_0+s_+)$ for suitable $s_l,s_+>0$ we see that $\tilde{\Gamma}_+$ is a smooth extension of $\Gamma|_{[s_0,s_0+s_+)}=:\Gamma_+$ which satisfies $\left|\dot{\tilde{\Gamma}}_+\right|_g=1$ on all of $\tilde{I}$ by construction. Similarly we can find a suitable $s_->0$ and a smooth extension $\tilde{\Gamma}_-$ of $\Gamma_-:=\Gamma|_{(s_0-s_-,s_0]}$ to some open interval $\tilde{I}$ around $s_0$ and which satisfies $\left|\dot{\tilde{\Gamma}}_-\right|_g=1$ at all points and $\tilde{\Gamma}_-(\tilde{I})\subseteq \gamma_0((0,\delta))\subseteq \psi((0,1))$. Since $L_n$ is smoothly embedded in $M$, we may view $\Gamma|_{(0,l)}$ as a map into $L_n$ and show that it is smooth. Then smoothness as a map into $M$ follows automatically. Thus take $\psi$ as a global chart of $L_n$ and consider the coordinate expressions $\tilde{\Gamma}_{\pm,loc}:=\psi^{-1}\circ \tilde{\Gamma}_{\pm}$. By choice of $\tilde{\Gamma}_{\pm}$ we know their tangent at $s_0$ is normalised and the corresponding tangent space is $1$-dimensional because these curves are contained in $L_n$. Assume for the moment that the tangents at $s_0$ point in opposite directions, i.e. $\dot{\tilde{\Gamma}}_{-,loc}(s_0)=-\dot{\tilde{\Gamma}}_{+,loc}(s_0)$. Then using a Taylor expansion around $\tau=0$ for $0<\tau$ small enough we find
\[
\tilde{\Gamma}_{-,loc}(s_0-\tau)=\tilde{\Gamma}_{-,loc}(s_0)-\dot{\tilde{\Gamma}}_{-,loc}(s_0)\tau+\ddot{\tilde{\Gamma}}_{-,loc}(\xi)\frac{\tau^2}{2}=\tilde{\Gamma}_{-,loc}(s_0)+\dot{\tilde{\Gamma}}_{+,loc}(s_0)\tau+\ddot{\tilde{\Gamma}}_{-,loc}(\xi)\frac{\tau^2}{2}
\]
for some $s_0-\tau\leq \xi\leq s_0$. Observe that $\tilde{\Gamma}_{+,loc}(s)=\sigma^{-1}(s)$ by (\ref{PT316}) and the choice of chart. Hence $\dot{\tilde{\Gamma}}_{+,loc}(s_0)>0$. In addition we take $\tau>0$ and the term $\ddot{\tilde{\Gamma}}_{-,loc}(\xi)$ can be uniformly bounded in $\tau$ for $\tau$ smaller than a fixed threshold. We conclude that for $0<\tau$ small enough we have $\dot{\tilde{\Gamma}}_{+,loc}(s_0)\tau+\ddot{\tilde{\Gamma}}_{-,loc}(\xi)\frac{\tau^2}{2}>0$ and since $\tilde{\Gamma}_-$ is an extension of the restriction of $\Gamma$ we compute $\tilde{\Gamma}_{-,loc}(s_0)=\epsilon_0$ and arrive at
\begin{equation}
\label{PT317}
\tilde{\Gamma}_{-,loc}(s_0-\tau)>\epsilon_0\text{ for all small enough }\tau>0.
\end{equation}
On the other hand since $\tilde{\Gamma}_-$ is an extension of $\Gamma_-$ and since $s_0-\tau<s_0$ we have $\tilde{\Gamma}_-(s_0-\tau)=\Gamma(s_0-\tau)=\hat{\gamma}(s^{-1}(s_0-\tau))$. Note that $s^{-1}$ is strictly increasing and $s^{-1}(s_0)=\epsilon_0$, thus $\tilde{\Gamma}_-(s_0-\tau)\in \gamma_0\left(\left(0,\frac{\delta}{2} \right)\right)=\psi((0,\epsilon_0))$ by definition of $\hat{\gamma}$ and $\epsilon_0$. Thus we find $0<\tilde{\Gamma}_{-,loc}(s_0-\tau)<\epsilon_0$, which contradicts (\ref{PT317}). We conclude that both tangent vectors at $s_0$ must point in the same direction and due to the normalisation condition already coincide. Since $\tilde{\Gamma}_{\pm}$ are extensions of the left and right restriction of $\Gamma$ respectively we see that $\Gamma$ is $C^1$ on $[0,l]$. In order to establish smoothness we will show that $\tilde{\Gamma}_+$ and $\tilde{\Gamma}_-$ coincide on some open neighbourhood around $s_0$, which will imply that $\Gamma$ coincides locally around $s_0$ with these smooth curves, i.e. is itself smooth around $s_0$. To see this we assume for the moment that there does not exist an open neighbourhood around $s_0$ such that $\tilde{\Gamma}_{\pm,loc}$ differ by only a constant on this neighbourhood. Then we can find a sequence $(\tau_n)_n$ converging to $s_0$ such that $\dot{\tilde{\Gamma}}_{+,loc}(\tau_n)\neq \dot{\tilde{\Gamma}}_{-,loc}(\tau_n)$ and due to the normalisation and one dimensionality of the tangent space we must have $\dot{\tilde{\Gamma}}_{+,loc}(\tau_n)=-\dot{\tilde{\Gamma}}_{-,loc}(\tau_n)$ for all $n$. Then a continuity argument in combination with the already established fact that $\dot{\tilde{\Gamma}}_{+,loc}(s_0)=\dot{\tilde{\Gamma}}_{-,loc}(s_0)$ implies $\dot{\tilde{\Gamma}}_{+,loc}(s_0)=-\dot{\tilde{\Gamma}}_{-,loc}(s_0)=-\dot{\tilde{\Gamma}}_{+,loc}(s_0)$, i.e. $\dot{\tilde{\Gamma}}_{+}(s_0)=0$, which contradicts the normalisation condition. Thus indeed in a small enough neighbourhood around $s_0$ both local expressions of the extensions differ by at most a constant. However we readily check that $\tilde{\Gamma}_{+,loc}(s_0)=\epsilon_0=\tilde{\Gamma}_{-,loc}(s_0)$ and hence this constant must be $0$. This shows that $\Gamma|_{(0,l)} \in C^{\infty}((0,l),M)$.
\newline
\newline
In order to establish the tameness of the knot $\text{clos}(L_n)$ it is left to show that $\Gamma$ has a finite total geodesic curvature. For notational simplicity we will simply write $\kappa$ instead of $\kappa_{\Gamma}$ for the geodesic curvature of $\Gamma$. We observe that since $\Gamma$ is smooth on the interior we have $\kappa\in C^0([a,b],\mathbb{R})$ for every $0<a<b<l$ and therefore it is enough to prove the existence of such $a,b$ which satisfy $\int_0^a\kappa ds,\int_b^l\kappa d\tau <\infty$. We will establish the existence of $a$, since the other case may be treated similarly. We recall the definition $\kappa(\tau)=\left|D_t\dot{\Gamma}(\tau) \right|_g$, i.e. $\kappa$ is the length of the acceleration vector of $\Gamma$, which in local coordinates can be expressed as
\[
D_t\dot{\Gamma}(\tau)=\left(\ddot{\Gamma}^i(\tau)+\dot{\Gamma}^j(\tau)\dot{\Gamma}^k(\tau)\Gamma^i_{jk}(\Gamma(\tau)) \right)\partial_i(\Gamma(\tau)),
\]
where $\Gamma^i_{jk}$ denotes the Christoffel symbols of the Levi-Civita connection. Now if we fix any chart $(\mu,U)$ around $p=\Gamma(0)$ with $\mu(p)=0$ and if we choose $a$ small enough, then $\Gamma([0,a])\subset U$. Now using the relation $2\alpha\beta\leq \alpha^2+\beta^2$, that $\Gamma\in C^1([0,l],M)$, that the metric is smooth and that $\Gamma([0,a])$ is compact, we can find a constant $c>0$ with
\begin{equation}
\label{PT318}
\kappa(\tau)\leq c\left(1+\sum_{i=1}^3|\ddot{\Gamma}^i(\tau)| \right)\text{ for all }0\leq \tau\leq a.
\end{equation}
Thus (\ref{PT318}) shows that it is enough to estimate $\int_0^a|\ddot{\Gamma}^i(\tau)|d\tau$ for $1\leq i \leq 3$. So fix any such $i$ and recall that $\Gamma=\hat{\gamma}\circ s^{-1}$. Since $s^{-1}$ is strictly monotonically increasing we may choose $a$ so small that $s^{-1}(\tau)< \epsilon_0$ for all $0\leq \tau \leq a$ so that we are in the first case of the definition of $\hat{\gamma}$, (\ref{PT315}), where we from now on set $\hat{\tilde{\gamma}}(t):=\tilde{\gamma}_0\left(\left(\frac{\delta}{2}\right)^{N_0}\frac{t}{\epsilon_0} \right)$. We recall that by construction $\tilde{\gamma}_0$ is $C^1$ up to the boundary and has a nowhere vanishing tangent. Then an explicit computation of $\ddot{\Gamma}^i$ by means of the chain rule allows us to estimate 
\begin{equation}
\label{PT319}
|\ddot{\Gamma}^i(\tau)|\leq C\left(1+\sum_{k=1}^3|\ddot{\hat{\tilde{\gamma}}}^k(s^{-1}(\tau))| \right),
\end{equation}
for a suitable $C>0$ independent of $\tau$. Hence it is enough to estimate $\int_0^a|\ddot{\hat{\tilde{\gamma}}}^k(s^{-1}(\tau))|d\tau$ for $1\leq k\leq 3$. To this end we perform a change of variables and can uniformly bound the Jacobian determinant since $\tilde{\gamma}_0$ is $C^1$ up to the boundary. We are eventually left with showing that $\int_0^{t_a}|\ddot{\tilde{\gamma}}_0^k(t)|dt<+\infty$ for suitably small $t_a$ or equivalently that $\int_0^{t_a}\sqrt{\sum_{k=1}^3|\ddot{\tilde{\gamma}}_0^k(t)|^2} dt$ is finite. To see this we have to be precise about how exactly the constant $N_0$ in the definition of (\ref{PT315}) was chosen, see also the defining equation (\ref{PT313}). We recall that $\gamma_0$ is a real analytic, non-constant curve defined on some open interval around $0$ with $\gamma_0(0)=p$. In our local coordinates it can be expressed for $0\leq |t|\ll 1$ as a convergent power series of the form
\begin{equation}
\label{PT320}
(\mu\circ \gamma_0)(t)=\sum_{k=1}^{\infty}a_kt^k, \text{ with }a_k\in \mathbb{R}^3
\end{equation}
and where we used $\mu(\gamma_0(0))=\mu(p)=0$. Since $\gamma_0$ is not constant there must be a smallest (strictly) positive integer $N_0$ with $a_{N_0}\neq 0$, which is exactly our choice of $N_0$. Setting $\vec{\tilde{\gamma}}_0:=(\tilde{\gamma}_0^1,\tilde{\gamma}_0^2,\tilde{\gamma}_0^3)$, we compute for $t>0$ and with the definition of $N_0$ in mind
\[
\ddot{\vec{\tilde{\gamma}}}_0(t)=\sum_{k=1}^{\infty}a_{N_0+k}\frac{k}{N_0}\left(\frac{k}{N_0}+1\right)t^{\frac{k}{N_0}-1}
\]
and consequently we can estimate by the triangle inequality
\begin{equation}
\label{PT321}
|\ddot{\vec{\tilde{\gamma}}}_0(t)|_2\leq \sum_{k=1}^{\infty}|a_{N_0+k}|_2\frac{k}{N_0}\left(1+\frac{k}{N_0}\right)t^{\frac{k}{N_0}-1}\text{ for all }0<t\ll 1.
\end{equation}
Since $\Sigma_n(t):=\chi_{(0,t_a)}(t)\sum_{k=1}^{n}|a_{N_0+k}|_2\frac{k}{N_0}\left(1+\frac{k}{N_0}\right)t^{\frac{k}{N_0}-1}$, where $\chi$ denotes the characteristic function, is a sequence of nonnegative, measurable functions which is monotonically increasing, we obtain by monotone convergence
\[
\int_0^{t_a}|\ddot{\vec{\tilde{\gamma}}}_0(t)|_2dt\leq \int_0^{t_a}\sum_{k=1}^{\infty}|a_{N_0+k}|_2\frac{k}{N_0}\left(1+\frac{k}{N_0}\right)t^{\frac{k}{N_0}-1}dt
\]
\[
=\sum_{k=1}^{\infty}\int_0^{t_a}|a_{N_0+k}|_2\frac{k}{N_0}\left(1+\frac{k}{N_0}\right)t^{\frac{k}{N_0}-1}dt=\sum_{k=1}^{\infty}|a_{N_0+k}|_2\left(1+\frac{k}{N_0}\right)t^{\frac{k}{N_0}}_a<+\infty
\]
for $t_a$ small enough. This shows the finiteness of the total geodesic curvature. Thus the last step is concluded and the proof of \cref{MT3} complete. $\square$
\newline
\newline
Keeping in mind the classification of $1$-manifolds, \cite{M65}, \cref{MC4} and \cref{MC5} are direct consequences of \cref{MT3}. $\square$
\section{Proof of \cref{MP7}}
In order to prove this result we start with a slight modification of the ABC-flows. To this end define $\mathcal{T}:=\mathbb{R}\slash (2\pi\mathbb{Z})\times \mathbb{R}\slash (2\pi\mathbb{Z}) \times \left(\frac{1}{2},\frac{3}{2}\right)$ and equip it with its standard differentiable structure and the flat metric $g_F$. This turns $\mathcal{T}$ into an orientable, real analytic $3$-dimensional manifold without boundary with a real analytic metric. We denote elements $P\in \mathcal{T}$ by $P=([x],[y],z)$, where $[\cdot]$ denotes the induced equivalence class. We then define
\begin{equation}
\label{PP71}
\psi: \mathcal{T}\rightarrow \mathcal{T}_A, P=([x],[y],z)\mapsto \left(\cos(x)(2+z\cos(y)),\sin(x)(2+z\cos(y)),z\sin(y) \right),
\end{equation}
where $\mathcal{T}_A$ is the solid toroidal annulus defined in (\ref{M2}). One readily checks that $\psi$ is well-defined and a real analytic diffeomorphism onto $\mathcal{T}_A$. Since $\mathcal{T}_A$ is an open subset of $\mathbb{R}^3$ it has a standard orientation and we choose the orientation on $\mathcal{T}$ such that $\psi$ is orientation preserving. We now define the following vector field $\bm{B}$ on $\mathcal{T}$ which is expressed in Cartesian coordinates by
\begin{equation}
\label{PP72}
\bm{B}(P):=\left(\cos(z-1)-\cos(y),\sin(1-z),-\sin(y) \right)\text{ for }P=([x],[y],z).
\end{equation}
One easily verifies that $\bm{B}$ is an eigenvector field of the curl operator on $(\mathcal{T},g_F)$, corresponding to the eigenvalue $+1$ with respect to the chosen orientation. Further one can confirm by direct calculations that
\begin{equation}
\label{PP73}
\mathcal{N}:=\{P\in \mathcal{T}|\bm{B}(P)=0 \}=\mathbb{R}\slash 2\pi\mathbb{Z}\times \{[0]\}\times \{1\}.
\end{equation}
Our construction proceeds as follows. Suppose we are given $(p,q)\in \mathbb{N}^2$, a pair of (strictly) positive coprime integers, then we denote by $b_0,d_0\in \mathbb{Z}$ the pair of integers generated by the extended Euclidean algorithm, which satisfies $pb_0+qd_0=1$. Given $k\in \mathbb{Z}$ we define $b_k:=b_0+kq$ and $d_k:=d_0-kp$ so that $pb_k+qd_k=1$ for all $k\in \mathbb{Z}$. We observe that the following map
\begin{equation}
\label{PP74}
f_{p,q,k}:\mathcal{T}\rightarrow \mathcal{T}, ([x],[y],z)\mapsto ([qx-b_ky],[px+d_ky],z)
\end{equation}
gives rise to a real analytic, orientation preserving diffeomorphism and we make the following definitions
\begin{equation}
\label{PP75}
\bm{X}_{p,q,k}:=(\psi\circ f_{p,q,k})_{*}\bm{B}\text{ and }g_{p,q,k}:=((\psi\circ f_{p,q,k})^{-1})^{\#}g_F,
\end{equation}
where $\cdot_{*}$ denotes the pushforward and $\cdot^{\#}$ denotes the pullback. Observe that all diffeomorphisms involved are orientation preserving and real analytic and that $\bm{B}$ is real analytic as well. Hence the above defined quantities are all real analytic and by construction of the $\bm{X}_{p,q,k}$ they are eigenfields of curl, with respect to $g_{p,q,k}$ and the standard orientation, corresponding to the eigenvalue $+1$. Further their zero set is by construction precisely given by
\[
\{(x,y,z)\in \mathcal{T}_A|\bm{X}_{p,q,k}((x,y,z))=0 \}=(\psi\circ f_{p,q,k})(\mathcal{N})=\mathcal{T}_{p,q},
\]
by (\ref{PP73}), \cref{MD6} and direct calculation. This proves properties (\ref{M3}) and (\ref{M4}). As for the explicit expressions one computes the following. Let
\begin{equation}
\label{PP76}
M_{p,q,k}:=\begin{pmatrix}
d^2_k+p^2 & d_kb_k-pq  & 0 \\
d_kb_k-pq & b^2_k+q^2 & 0 \\ 
0 & 0 & 1
\end{pmatrix}\text{ and } D(x,y,z):=\begin{pmatrix}
-\frac{y}{r^2} & \frac{x}{r^2}& 0\\
-\frac{xz}{rR^2} & -\frac{zy}{rR^2} & \frac{(r-2)}{R^2}\\
\frac{x(r-2)}{rR}& \frac{y(r-2)}{rR} & \frac{z}{R}
\end{pmatrix},
\end{equation}
where $r:=\sqrt{x^2+y^2}$ and $R:=\sqrt{(r-2)^2+z^2}$. Then $g_{p,q,k}$ is represented in Euclidean coordinates by the matrix
\begin{equation}
\label{PP77}
g_{p,q,k}(x,y,z)=D^{Tr}(x,y,z)\cdot M_{p,q,k}\cdot D(x,y,z),
\end{equation}
where $\cdot$ indicates the standard matrix multiplication. It follows easily from (\ref{PP77}) that the metrics are indeed distinct for different values of $k$ since the matrix $D(x,y,z)$ is invertible for all points in $\mathcal{T}_A$. As for $\bm{X}_{p,q,k}$ we have the following: Given $(x,y,z)\in \mathcal{T}_A$ we define $P=([a],[c],t):=\psi^{-1}((x,y,z))$, then we have the expression
\[
\bm{X}_{p,q,k}(x,y,z)=\left(qX_{p,q}^1(P)-b_kX_{p,q}^2(P)\right)\begin{pmatrix}
-y \\
x & \\
0& 
\end{pmatrix}
\]
\begin{equation}
\label{PP78}
+\left(pX_{p,q}^1(P)+X_{p,q}^2(P)d_k\right)\begin{pmatrix}
 -\frac{zx}{r}\\
 -\frac{zy}{r} \\
 r-2 
\end{pmatrix}+X_{p,q}^3(P)\begin{pmatrix}
\frac{x(r-2)}{rR}\\
\frac{y(r-2)}{rR}\\
\frac{z}{R}
\end{pmatrix},
\end{equation}
where $r=r(x,y,z)$ and $R=R(x,y,z)$ are defined as before and where
\begin{equation}
\label{PP79}
\begin{pmatrix}
X_{p,q}^1(P)\\
X_{p,q}^2(P)\\
X_{p,q}^3(P)
\end{pmatrix}:=\begin{pmatrix}
\cos(R(x,y,z)-1)-\left(\cos(qc)\cos(pa)+\sin(qc)\sin(pa) \right)\\
\sin(1-R(x,y,z))\\
\cos(qc)\sin(pa)-\sin(qc)\cos(pa)
\end{pmatrix}.
\end{equation}
Observe that $p$ and $q$ are positive integers and hence by means of standard trigonometric identities it is always possible to express the quantities $\cos(qc),\cos(pa),\sin(qc),\sin(pa)$ as polynomials in the variables $\cos(a),\sin(a),\cos(c),\sin(c)$ for which we have the following identities in terms of Cartesian coordinates
\begin{equation}
\label{PP710}
\cos(c)=\frac{r-2}{R},\text{ }\sin(c)=\frac{z}{R},\text{ }\cos(a)=\frac{x}{r},\text{ }\sin(a)=\frac{y}{r},
\end{equation}
with the usual functions $r=r(x,y,z)$ and $R=R(x,y,z)$. Thus we obtain from (\ref{PP78})-(\ref{PP710}) an explicit expression for the $\bm{X}_{p,q,k}$ in terms of Cartesian coordinates, which involves only elementary functions. Lastly note that by properties of the pushforward and since $\psi$ is a diffeomorphism, the constructed vector fields $\bm{X}_{p,q,k}$ are all distinct if and only if the vector fields $(f_{p,q,k})_{*}B$ are all distinct, which is easy to verify. This concludes the proof of \cref{MP7}. $\square$
\begin{exmp}[The trefoil knot]
	\label{PP7E1}
	The simplest, non-trivial torus knot is the $(p,q)=(2,3)$ torus knot, also known as the trefoil knot. From the Euclidean algorithm we get $b_0=-1$, $d_0=1$ and by means of trigonometric identities, as described above, we can compute the quantities $X^j_{2,3}$ in (\ref{PP79}), omitting the arguments, to be
	\begin{equation*}
	X^1_{2,3}=\cos(R-1)+3\frac{x^2z^2(r-2)}{r^2R^3}+\frac{r-2}{R}\left(\frac{y^2}{r^2}\left(2-5\frac{z^2}{R^2}\right)-\frac{(r-2)^2}{R^2} \right)-2\frac{xyz}{r^2R}\left(4\frac{(r-2)^2}{R^2}-1\right)
	\end{equation*}
	\begin{equation*}
	X^2_{2,3}=\sin(1-R)
	\end{equation*}
	\begin{equation*}
	X^3_{2,3}=\frac{x^2z^3}{r^2R^3}+\frac{z}{r}\left(\frac{y^2}{r^2}\left(7\frac{(r-2)^2}{R^2}-1 \right)-3\frac{(r-2)^2}{R^2} \right)+2\frac{xy(r-2)}{r^2R}\left(1-4\frac{z^2}{R^2} \right),
	\end{equation*}
	where as usual $r=\sqrt{x^2+y^2}$ and $R=\sqrt{(r-2)^2+z^2}$. By definition we have $b_k=3k-1$ and $d_k=1-2k$, so that we obtain an exact expression for the family of vector fields $\bm{X}_{2,3,k}$ all of whose zero sets are given by the same trefoil knot. \Cref{Figure1} depicts the corresponding zero set as a subset of the solid toroidal annulus $\mathcal{T}_A$, which is a solid torus from whose interior a smaller solid torus was cut out.
\begin{figure}[htbp] 
	\centering
	\includegraphics[width=0.5\textwidth]{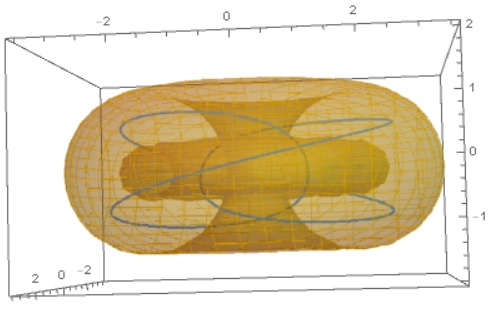}
	\caption{The zero set of the vector fields $\bm{X}_{2,3,k}$ in blue lying within the solid toroidal annulus, depicted in yellow}
	\label{Figure1}
\end{figure}
\end{exmp}
\section*{Acknowledgements}
This work has been funded by the Deutsche Forschungsgemeinschaft (DFG, German Research Foundation) – Projektnummer 320021702/GRK2326 –  Energy, Entropy, and Dissipative Dynamics (EDDy). I would further like to thank Christof Melcher and Heiko von der Mosel for discussions.
\bibliographystyle{plain}
\bibliography{mybibfile}

\begin{thebibliography}{10}

\bibitem{A66}
V.I. Arnold.
\newblock \foreignlanguage{russian}{О топологии трехмерных
  стационарных течений идеальной жидкости}
  (on the topology of three-dimensional steady flows of an ideal fluid).
\newblock {\em PMM}, 30(1):183--185, 1966.

\bibitem{A74}
V.I. Arnold.
\newblock The asymptotic {H}opf invariant and its applications.
\newblock In {\em Proceedings of the All-Union Summer School on Partial
  Differential Equations (Dilizhan, Erevan, Armenia)}, pages 229--256. Armenian
  SSR Academy of Sciences Press, 1974.

\bibitem{AK98}
V.I. Arnold and B.A. Khesin.
\newblock {\em Topological Methods in Hydrodynamics.}
\newblock Springer Verlag, 1998.

\bibitem{AL91}
M.~Avellaneda and P.~Laurence.
\newblock On {W}oltjer's variational principle for force-free fields.
\newblock {\em Journal of Mathematical Physics}, 32(5):1240--1253, 1991.

\bibitem{BM88}
E.~Bierstone and P.D. Milman.
\newblock Semianalytic and subanalytic sets.
\newblock {\em Publications Math\'{e}matiques de l'IH\'{E}S}, 67:5--42, 1988.

\bibitem{C99}
J.~Cantarella.
\newblock {\em Topological structure of stable plasma flows}.
\newblock PhD thesis, University of Pennsylvania, 1999.

\bibitem{DFGHMS86}
T.~Dombre, U.~Frisch, J.M. Greene, M.~H\'{e}non, A.~Mehr, and A.M. Soward.
\newblock Chaotic streamlines in the {A}{B}{C} flows.
\newblock {\em Journal of Fluid Mechanics}, 167:353--391, 1986.

\bibitem{EP12}
A.~Enciso and D.~Peralta-Salas.
\newblock Knots and links in steady solutions of the {E}uler equation.
\newblock {\em Ann. of Math.}, 175(1):345--367, 2012.

\bibitem{EP15}
A.~Enciso and D.~Peralta-Salas.
\newblock Existence of knotted vortex tubes in steady {E}uler flows.
\newblock {\em Acta Math.}, 214:61--134, 2015.

\bibitem{EP16}
A.~Enciso and D.~Peralta-Salas.
\newblock {B}eltrami fields with a nonconstant proportionality factor are rare.
\newblock {\em Arch. Rational Mech. Anal.}, 220:243--260, 2016.

\bibitem{EtGr00a}
J.~Etnyre and R.~Ghrist.
\newblock Contact topology and hydrodynamics i: {B}eltrami fields and the
  {S}eifert conjecture.
\newblock {\em Nonlinearity}, 13(2):441--458, 2000.

\bibitem{EtGr00b}
J.~Etnyre and R.~Ghrist.
\newblock Contact topology and hydrodynamics iii: Knotted flowlines.
\newblock {\em Transactions of the American Mathematical Society},
  352(12):5781--5794, 2000.

\bibitem{G68}
A.M. Gabri\'{e}lov.
\newblock \foreignlanguage{russian}{О проекциях
  полуаналитических множеств} (projections of
  semi-analytic sets).
\newblock {\em Funktsional'nyi Analiz i Ego Prilozheniya}, 2(4):18--30, 1968.

\bibitem{HarSi89}
R.~Hardt and L.~Simon.
\newblock Nodal sets for solutions of elliptic equations.
\newblock {\em J. differential Geometry}, 30(2):505--522, 1989.

\bibitem{H093}
H.~Hofer.
\newblock Pseudoholomorphic curves in symplectizations with applications to the
  {W}einstein conjecture in dimension three.
\newblock {\em Inventiones mathematicae}, 114:515--563, 1993.

\bibitem{KrPa02}
S.~Krantz and H.~Parks.
\newblock {\em A Primer of Real Analytic Functions}.
\newblock Birkh{\"a}user Basel, second edition, 2002.

\bibitem{Log18}
A.~Logunov.
\newblock Nodal sets of {L}aplace eigenfunctions: polynomial upper estimates of
  the {H}ausdorff measure.
\newblock {\em Annals of Mathematics}, 187(1):221--239, 2018.

\bibitem{L65}
S.~\L{}ojasiewicz.
\newblock {\em Ensembles semi-analytiques}.
\newblock Institut des Hautes Etudes Scientifiques, Bures-sur-Yvette, 1965.

\bibitem{Loo84}
E.~Looijenga.
\newblock {\em Isolated Singular Points on Complete Intersections}.
\newblock Cambridge University Press, 1984.

\bibitem{MaTr07}
D.~Massey and L.D. Tr\'{a}ng.
\newblock Notes on real and complex analytic and semianalytic singularities.
\newblock In {\em Singularities in Geometry and Topology Proceedings of the
  Trieste Singularity Summer School and Workshop}, pages 81--126. World
  Scientific Publishing, 2007.

\bibitem{M50}
J.W. Milnor.
\newblock On the total curvature of knots.
\newblock {\em The Annals of Mathematics}, 52(2):248--257, 1950.

\bibitem{M65}
J.W. Milnor.
\newblock {\em Topology from the differentiable view point}.
\newblock Princeton University Press, 1965.

\bibitem{M68}
J.W. Milnor.
\newblock {\em Singular Points of Complex Hypersurfaces}.
\newblock Princeton University Press, 1968.

\bibitem{N14}
N.~Nadirashvili.
\newblock {L}iouville theorem for {B}eltrami flow.
\newblock {\em Geom. Funct. Anal.}, 24:916--921, 2014.

\bibitem{SoZel11}
C.D. Sogge and S.~Zelditch.
\newblock Lower bounds on the {H}ausdorff measure of nodal sets.
\newblock {\em Math. Res. Lett.}, 18(1):25--37, 2011.

\bibitem{SoZel12}
C.D. Sogge and S.~Zelditch.
\newblock Lower bounds on the {H}ausdorff measure of nodal sets ii.
\newblock {\em Math. Res. Lett.}, 19(6):1361--1364, 2012.

\bibitem{Sul08}
J.M. Sullivan.
\newblock Curves of finite total curvature.
\newblock In {\em Discrete Differential Geometry}, pages 137--161.
  Birkh\"{a}user Basel, 2008.

\bibitem{Ta07}
C.H. Taubes.
\newblock The {S}eiberg-{W}itten equations and the {W}einstein conjecture.
\newblock {\em Geom. Topol.}, 11:2117--2202, 2007.

\bibitem{W58}
L.~Woltjer.
\newblock A theorem on force-free magnetic fields.
\newblock In {\em Proc. Natl. Acad. Sci. USA}, volume~44, pages 489--491, 1958.

\end{thebibliography}
\footnotesize
\end{document}